\documentclass[]{article}
\usepackage{amsfonts}
\usepackage{theorem}
\newtheorem{construction}{Construction}

\newtheorem{theorem}{Theorem}

\newtheorem{lemma}{Lemma}
\newtheorem{remark}{Remark}
\newcommand{\openbox}{$\begin{array}{c}
\hspace*{-0.55em}\sqcap \hspace*{-0.60em}\\[-0.4em] \hline
\multicolumn{1}{c}{\hspace*{-0.60em}}\\[-0.8em]
\end{array}
$}

\usepackage{enumerate}

\begin{document}

\centerline{\bf Finite permutable Putcha semigroups
\footnote{
The published form of this paper is available at\\
http://acta.fyx.hu/acta/showCustomerArticle.action?id=13007\&dataObjectType=\\
article\&returnAction=showCustomerVolume\&sessionDataSetId=7fcdd07e30c28c3\&style=
}}

\medskip

\centerline{A. De\'ak and A. Nagy}

\bigskip

%\centerline{Communicated by M. B. Szendrei}

\bigskip

\begin{abstract}A semigroup $S$ is called a permutable semigroup if
$\alpha \circ \beta =\beta \circ \alpha$ is satified for all
congruences $\alpha$ and $\beta$ of $S$. A semigroup is called a Putcha semigroup
if it is a semilattice of archimedean semigroups. In this paper we deal with
finite permutable Putcha semigroups. We describe the finite permutable
archimedean semigroups and finite permutable semigroups which are
semilattices of a group and a nilpotent semigroup.

\medskip

\end{abstract}

\section{Introduction}

The notion of the permutable semigroup was introduced in \cite{Hamilton}. A semigroup $S$ is called
a permutable semigroup if $\alpha \circ \beta =\beta \circ \alpha$ is satisfied for all congruences
$\alpha$ and $\beta$ of $S$. In \cite{Hamilton}, H. Hamilton proved some general results about
permutable semigroups and described the commutative permutable semigroups.
In \cite{Cherubini medial}, \cite{Bonzini typ2}, \cite{Bonzini typ1}, \cite{Cherubini duo},
\cite{Cherubini reg}, \cite{Cherubini complreg}, \cite{Deak}, \cite{Jiang},
\cite{Jiang-Chen}, \cite{Nagy non-trivial}, \cite{Nagy typ1}
the permutable semigroups are examined in special
classes of semigroups. In the present paper we deal with finite permutable Putcha semigroups, that is,
finite semigroups which are semilattices of archimedean semigroups (a semigroup $S$
is called archimedean if, for every $a,b\in S$, there are positive integers $n$ and $m$
such that $a^n\in SbS$ and $b^m\in SaS$).
We show that the finite archimedean permutable semigroups are exactly the finite cyclic
nilpotent semigroups and the finite completely simple
permutable semigroups. Dealing with the non-archimedean case, we  describe
such finite permutable semigroups which are semilattices of a group $G$ and a nilpotent
semigroup $N$ with $GN\subseteq N$.

First of all we cite some earlier results which will be used in our investigation.
For notations and notions not defined here, we refer to \cite{Clifford} and \cite {Nagy book}.

\begin{lemma}\label{L1} Every group is permutable.\hfill\openbox
\end{lemma}

\begin{lemma}\label{L2}(\cite {Hamilton}) The ideals (equivalently, the principal ideals) of a permutable
semigroup form a chain with respect to inclusion.\hfill\openbox
\end{lemma}

\begin{lemma}\label{L3}(\cite{Nagy non-trivial}) A nil semigroup $S$ is permutable if and only if the ideals of $S$
form a chain with respect to inclusion.\hfill\openbox
\end{lemma}

\begin{lemma}\label{L4}(\cite{Hamilton}) If a permutable semigroup $S$ contains a proper ideal $K$ then neither $S$
nor $K$ has a non-trivial group homomorphic image.\hfill\openbox
\end{lemma}

\begin{lemma}\label{L5}(\cite{Jiang}) If $\rho$ is an arbitrary congruence and $K$ is an arbitrary ideal
of a permutable semigroup then $K$ is in an $\alpha$-class or $K$ is a union
of $\alpha$-classes.\hfill\openbox
\end{lemma}

\begin{lemma}\label{L6}(\cite{Cherubini complreg}) A Rees matrix semigroup $S={\cal M}(G,I,J;P)$ is permutable if
and only if $|I|\leq 2$ and $|J|\leq 2$.\hfill\openbox
\end{lemma}

\begin{lemma}\label{L7}(\cite{Hamilton}) Every homomorphic image of a permutable semigroup is permutable.\hfill\openbox
\end{lemma}

\begin{lemma}\label{L8}(\cite{Hamilton}) A semilattice $F$ is permutable if and only if $|F|\leq 2$.\hfill\openbox
\end{lemma}

\medskip

As a Putcha semigroup is a semilattice of archimedean semigroups, Lemma \ref{L7}
and Lemma \ref{L8} together
imply that a permutable Putcha semigroup is either archimedean or a semilattice of two
archimedean semigroups $S_0$ and $S_1$ such that, for example, $S_0S_1\subseteq S_0$.

\section{Finite permutable archimedean semigroups}

First of all we formulate some lemmas which will be used in our examinations several times.

A semigroup $S$ with a zero $0$ is called a nil semigroup if, for every $a\in S$, there is
a positive integer $n$ such that $a^n=0$. $S$ is called nilpotent if, $S^m=\{0\}$ for some
positive integer $m$.

It is clear, that every archimedean semigroup with a zero is a nil semigroup.

\begin{lemma}\label{L9} Every finite nil semigroup is nilpotent.
\hfill\openbox
\end{lemma}

It is obvious, that every finite semigroup has a kernel $K$ which is completely simple.
Thus we have the following lemma.

\begin{lemma}\label{L10} A finite semigroup is archimedean if and only if it
is an ideal extension of a completely
simple semigroup by a nilpotent semigroup.\hfill\openbox
\end{lemma}

\begin{theorem}\label{T1} A finite semigroup is an archimedean permutable semigroup if and only if it is
either a cyclic nilpotent semigroup or a permutable completely simple semigroup.
\end{theorem}

\noindent
{\bf Proof}. Let $S$ be a finite permutable archimedean semigroup. By Lemma \ref{L10},
$S$ is an ideal extension of a completely simple semigroup $K$ ($K$ is the
kernel of $S$) by the nilpotent semigroup $N=S/K$. By Lemma \ref{L7} of this paper, Lemma 2 of \cite{Nagy non-trivial}
and Lemma 2 of \cite{Nagy Jones},
$N$ is a cyclic nilpotent semigroup. If $|K|=1$ then $S$ is isomorphic to $N$.
Consider the case when $|K|>1$. We show that $S=K$. Assume, in an indirect way, that $K\neq S$.
As $K$ is completely simple, the Green's
relations $\cal R$ and $\cal L$ are congruences on $K$. It is easy to see that
${\cal R}\cup 1_S$ and  ${\cal L}\cup 1_S$ are congruences on $S$ and the kernels
of the quotient semigroups are, respectively, the left zero semigroup $K/{\cal R}$
and the right zero semigroup $K/{\cal L}$. By Lemma \ref{L4}, $K/{\cal R}$ or
$K/{\cal L}$ is non-trivial. By symmetry, it can be assumed without loss of generality
that $K$ is a non-trivial right zero semigroup. Let $a\in S-K$ and let
$f=a^n\in K$, so $fa^i=f$ for all positive integers $i$ and $xf=f$ for all $x\in S$. Let $b\in K$, $b\neq f$.
Applying Lemma \ref{L5}, $b$ is related to $f$ under the congruence $\rho$ on $S$
generated by $(a,f)$, so there exists a sequence of elementary $\rho$-transitions from
$b$ to $f$ that begins either $b=sat\ \mapsto \ sft$ or $b=sft\ \mapsto \ sat$
($s,t\in S^1$), where the right hand side is distinct from $b$. In addition, since $b=bb$ and
$f=bf$, we can assume without loss of generality that $s=bs\in K$. If $t=1$ then $b=sa$
(since $b\neq bf=f$); otherwise, since $K$ is right zero, $t\notin K$, so $t=a^i$ for
some $i<n$ and $b=sa^{i+1}$, since again $b\neq sfa^i=f$. In either case, $b=ca$ for some
$c\in K$, $c\neq f$. Now the same argument applies to $c$ and iterating the argument leads
to $b=xa^n=xf=f$, a contradiction. Thus the first part of the theorem is proved.
As the converse is obvious, the theorem is proved.\hfill\openbox

\section{Finite permutable non-archimedean Putcha \newline semigroups}

\begin{lemma}\label{L11} If $S$ is a finite non-archimedean Putcha permutable semigroup then
it is a semilattice of a completely simple semigroup $S_1=M(G;I,J;P)$ such that $|I|,|J|\leq 2$
and a semigroup $S_0$ such that $S_1S_0\subseteq S_0$ and $S_0$ is an ideal extension
of a completely simple semigroup $K$ by a nilpotent semigroup.
\end{lemma}

\noindent
{\bf Proof}. Let $S$ be a finite permutable non-archimedean Putcha-semigroup. Then, by Lemma \ref{L7}
and Lemma \ref{L8},
$S$ is a semilattice of two archimedean semigroups $S_0$ and $S_1$ such that
$S_0S_1\subseteq S_0$. As the Rees factor $S_1^0=S/S_0$ is permutable by Lemma \ref{L7},
$S_1$ is a permutable archimedean semigroup. By Lemma 8 of \cite{Hamilton} and
Theorem 1 of this paper, $S_1$ is completely simple.
Then $S_1$ is a Rees matrix semigroup $S_1=
(G;I,J;P)$ and $|I|, |J|\leq 2$ by Lemma \ref{L6}.
By Lemma \ref{L10}, $S_0$ is an ideal extension of a completely simple semigroup $K$
by the nilpotent Rees factor semigroup $N=S_0/K$.\hfill\openbox

\medskip

In this paper we deal with only that case when $S_1$ is a group.

\begin{lemma}\label{L12} If a finite permutable semigroup $S$ is a semilattice of a group $G$
and a nilpotent semigroup $N$ such that $NG\subseteq N$ then the identity element of $G$
is a left identity element or a right identity element of $S$.
\end{lemma}

\noindent
{\bf Proof}.
Let $a\in N$ be an arbitrary element.
Then $J(a)\subseteq J(e)$, where $e$ denotes
the identity element of $G$.
Then there are elements $x,y\in S^1$ such that
$a=xey$. So $N=eN\cup Ne\cup NeN$. Since $N$ is an ideal, $Ne\cup NeN$ and $eN\cup NeN$
are ideals of $S$ and so, by hypothesis, one is included in the other. Suppose
$eN\subseteq Ne\cup NeN$, so that
$$N=Ne\cup NeN=Ne\cup (Ne)(eN)\subseteq Ne\cup (Ne)(Ne\cup NeN)\subseteq Ne\cup (Ne)^2N.$$
Inductively, $N\subseteq Ne\cup (Ne)^iN$ for all positive integers $i$, and since
$N$ is nilpotent, $N=Ne$, as required.
In case $Ne\subseteq eN\cup NeN$, we get $N=eN$.\hfill\openbox

\begin{lemma}\label{L13} Let $S$ be a finite permutable semigroup which is a
semilattice of a group $G$ and a
nilpotent semigroup $N$ such that $GN\subseteq N$. Let $e$ denote the identity
element of $G$. If $Ne=N$ then
$eN=\{ 0\}$ or $eN=N$. Similarly, if $eN=N$ then $Ne=\{ 0\}$ or $Ne=N$.
\end{lemma}

\noindent
{\bf Proof}: By the symmetry, we deal with only the first assertion of the lemma.
Assume $N=Ne$. Then $SeN=eN\cup NeN$, which is an ideal of $S$.
If $SeN=N$, then $N=eN\cup NeN$ from which we get $N=eN$ as in the proof of Lemma \ref{L12}.
If $SeN\neq N$, then consider the equivalence
$$\alpha=\{ (a,b)\in S\times S:\ ea=eb\}.$$
It is obvious that $\alpha$ is a right congruence.
Let $a,b,s$ be arbitrary elements of $S$
such that $(a,b)\in \alpha$. As $e$ is a right identity element of $S$, we get
$$sa=(se)a=s(ea)=s(eb)=(se)b=sb$$ and so
$$e(sa)=e(sb).$$ Thus $\alpha$ is also a left congruence of $S$, and so it is a congruence
of $S$. Let $x\in N$ be an arbitrary element. As $(x,ex)\in \alpha$ and $ex\in SeN$,
by Lemma \ref{L5}, the ideal $SeN$ is contained by the $\alpha$-class of $x$,
and so $(0,x)\in \alpha$, that is, $0=e0=ex$. Hence $eN=\{ 0\}$.
Thus the lemma is proved.\hfill\openbox

\begin{lemma}\label{L14} Let $S$ be a finite non-archimedean permutable semigroup which is
a semilattice of a group $G$ and an archimedean semigroup $S_0$ such that $GS_0\subseteq S_0$.
Then $S_0$ is either
\begin{itemize}
\item[(1)] completely simple,
\item[(2)] or a non-trivial null semigroup $N$ such that the identity element of
$G$ is a right identity element of $S$ and $SN=\{ 0\}$,
\item[(3)] or a non-trivial null semigroup $N$ such that the identity element of
$G$ is a left identity element of $S$ and $NS=\{ 0\}$,
\item[(4)] or an ideal extension of
a completely simple semigroup $K$ by a non-trivial nilpotent semigroup $N$ such that
the identity element of $G$ is an identity element of the factor semigroup $S/K$.
\end{itemize}
\end{lemma}

\noindent
{\bf Proof}. By Lemma \ref{L10}, $S_0$ is an ideal extension of a completely simple semigroup $K$ by the
Rees factor semigroup $N=S_0/K$ which is nilpotent. If $S_0=K$ then (1) is satisfied.

Assume $S_0\neq K$. As $K=K^2$ is an ideal of $S_0$ and $S_0$ is an ideal of $S$, we have that
$K$ is an ideal of $S$ (see Exercise 4(a) for \S 2.9 of \cite{Clifford}). Consider the Rees factor semigroup $S/K$ which
is a semilattice of $G$ and a nilpotent semigroup which is isomorphic
to the non-trivial semigroup $N=S_0/K$. By Lemma \ref{L12}, the identity element of $G$ is the right identity element or the left identity
element of $S/K$.

First consider the case when the identity element $e$ of $G$ is the right identity
element but not a left identity element of $S/K$. Then, by Lemma \ref{L13}, $eS_0\subseteq K$.
Now without loss of generality, if $K$ is non-trivial it can be assumed to be
either left zero or right zero, but the two cases must be treated separately
because of the asymmetry of the hypothesis on $S$. In either case, let
$a\in S-K-G$, such that $f=a^2\in K$, and suppose $b\in K$, $b\neq f$.
By Lemma \ref{L5}, $b$ is related to $f$ under the congruence $\rho$ on $S$
generated by $(a,f)$, so there exists a sequence of elementary $\rho$-transitions from
$b$ to $f$ that begins either $b=sat\ \mapsto \ sft$ or $b=sft\ \mapsto \ sat$
($s,t\in S^1$), where the right hand side is distinct from $b$.
First suppose that $K$ is right zero. Then again $t\notin K$. If $t\in N$ then
$at\in K$ and so $at=a(at)=ft$, giving $sat=sft$, a contradiction. So $t\in G$ and therefore
$b=be$. Hence $K=Ke$. As in the proof of Theorem \ref{T1}, without loss of generality, $s\in K$ and so $s=se$.
Also $ea\in K$. Then $$sa=(se)a=s(ea)=ea=a(ea)=(ae)a=a^2=f=sf,$$ again giving the
contradiction $sat=sft$.
Next suppose $K$ is left zero. Now, without loss of generality, $t\in K$ and $s\notin K$.
If $s\neq 1$ then since $S\subseteq K$, $sa=sasa=sa^2=sf$ (since $sas=sa$), a
contradiction. So $s=1$ and since $b\neq f=ft$, $b=at$. But $t\in K$ and $t\neq f$
(since $af=f$) so similarly $t=at'$ for some $t'$, yielding the contradiction
$b=a^2t'=ft'=f$. From this it follows, that $|K|=1$ and so $S_0=N$.
Let $a\in S_0=N$ be arbitrary.
As $eN=\{ 0\}$ and $ea\in eN$, we get $ea=0$ and so, for every $s\in S$, $sa=sea=0$.
Thus $SN=\{ 0\}$ and so (2) is satisfied.

If the identity element of $G$ is a left identity element but not a right identity element
of $S/K$ then (3), the dual of (2) is satisfied.

If the identity element of $G$ is the identity element of $S/K$, then (4) is satisfied.
Thus the lemma is proved.\hfill\openbox

\begin{remark} \label{R1}
If $|S_0|=1$ is satisfied in case $(1)$ of Lemma \ref{L14} then $S$ is a group with a zero adjoined
and so $S$ is permutable.
\end{remark}

\begin{remark}\label{R2} Condition $(4)$ of Lemma \ref{L14} has two subcases:

(4a): $|K|=1$ and so $S_0$ is a non-trivial nilpotent semigroup such that
the identity element of $G$ is an identity element of $S$.

(4b): $|K|>1$, but $K\neq S_0$.
\end{remark}

In this paper we describe only those finite permutable non-archimedean
Putcha semigroups which are
semilattice of a group $G$ and a semigroup $S_0$ with $GS_0\subseteq S_0$, where
$S_0$ satisfies either
condition $(2)$ or condition $(3)$ of Lemma \ref{L14} or condition $(4a)$ of
Remark \ref{R2}.

\subsection{When the identity element of $G$ is only a one-sided identity
element of $S$}

In this section we deal with only the right side case, but the main theorem (Theorem \ref{T2})
will be formulated for both right and left cases.

\medskip

\noindent
For a non-trivial nil semigroup $N$, let $N^*$ denote $N-\{ 0\}$.

\begin{lemma}\label{L15} Let $S$ be a finite permutable semigroup which is a semilattice of
a group $G$ and a non-trivial nilpotent semigroup $N$ such that the identity element of $G$ is
a right identity element of $S$ and $SN=\{ 0\}$. Then $aG=N^*$ for every
$a\in N^*$.
\end{lemma}

\noindent
{\bf Proof}.
It is clear that $N^2=\{ 0\}$. Let $a\in N^*$ be
arbitrary. If $ag=0$ for some $g\in G$ then $a=ae=agg^{-1}=0$ which is a contradiction. Thus
$aG\subseteq N^*$. As $S$ is finite and the ideals of $S$ form a chain, there
is an element $b\in N^*$ such that $S^1bS^1=N$. Thus
$N=S^1bS^1=(S\cup 1)b(G\cup N\cup 1)=bG\cup bN=bG\cup \{ 0\}$. Thus $bG=N^*$.
From this it follows that, for an arbitrary $a\in N^*$ and some $x\in G$,
$aG=bxG=bG=N^*$. \hfill\openbox

\begin{remark}\label{R3}
If a semigroup $S$ satisfies the conditions of Lemma \ref{L15} then $N^*$ is
a right $G$-set (\cite{McKenzie}) and $G$ acts on $N^*$ transitively.
\end{remark}

\begin{lemma}\label{L16} If an arbitrary semigroup $S$ is a semilattice of a group $G$ and a non-trivial
null semigroup $N$ such that $GN=\{ 0\}$ and $aG=N^*$ for every $a\in N^*$ then,
for every non-universal congruence $\alpha$ of $S$,
$[0]_{\alpha}$ is either $\{ 0\}$ or $N$, and
$[g]_{\alpha}\subseteq G$ for every $g\in G$.
\end{lemma}

\noindent
{\bf Proof}.
If $g\in [0]_{\alpha}$ for some $g\in G$ then $G\subseteq [0]_{\alpha}$ and so
$N^*=aG\subseteq [0]_{\alpha}$, where $a\in N^*$ is an arbitrary element.
Then $[0]_{\alpha}=S$. If $\alpha$ is not a universal congruence of $S$ then
$[0]_{\alpha}\subseteq N$. Assume $[0]_{\alpha}\neq \{ 0\}$. Then
there is an element $a\in N^*$ such that $a\in [0]_{\alpha}$, and so
$N^*=aG\subseteq I$. Hence $I=N$.

Assume $(a,g)\in \alpha$ for some $a\in N, g\in G$ and a non-universal
congruence $\alpha$ of $S$. Then $(ea,g)\in \alpha$, where $e$ is the
identity element of $G$. As $ea=0$, we get $(0,g)\in \alpha$ and so
$(0,h)\in \alpha$ for every $h\in G$ and  so $\alpha$ is the universal congruence
of $S$ by the above. It is a contradiction. Thus $[g]_{\alpha}\subseteq G$ for
every $g\in G$.\hfill\openbox

\begin{remark}\label{R4} By Lemma 4.20 of \cite{McKenzie}, if $X$ is a right $G$-set such that the group $G$ acts
on the non-empty set $X$ transitively then the congruence lattice $Con(X)$ of
the $G$-set $X$ is isomorphic to the interval $[Stab_G(x),G]$ for every $x\in X$,
where $Stab_G(x)=\{ g\in G:\ xg=x\}$. The corresponding isomorphisms are
$$\phi :\alpha \ \mapsto H_{\alpha}=\{ g\in G:\ xg\ \alpha \ x\}\ (\alpha \in {\it Con(X)})$$
and
$$\psi : H\ \mapsto \ \alpha _H=\{ (xg,xh)\in X\times X:\ Hg=Hh\}\ (H\in [Stab_G(x),G])$$
(which are inverses of each other).
\end{remark}

\medskip

As two right congruences of a group $G$ determined by
subgroups $H$ and $K$ of $G$ commute with each other if and only if $HK=KH$
it is easy to see (by Remark \ref{R4}) that the following lemma is true.

\begin{lemma}\label{L17} Let $X$ be a right $G$-set such that $G$ acts on $X$ transitively.
Let $x\in X$ be an arbitrary fixed element. Then $\alpha \circ \beta=
\beta \circ \alpha$ is satisfied for some congruences $\alpha, \beta \in Con(X)$
if and only if
$H_{\alpha}H_{\beta}=H_{\beta}H_{\alpha}$ is satisfied for $H_{\alpha}, H_{\beta}\in [Stab_G(x),G]$.
\hfill\openbox
\end{lemma}

\begin{remark}\label{R5}
As the congruence lattice $Con(X)$ of
a $G$-set $X$ is isomorphic to the interval $[Stab_G(x),G]$ for all $x\in X$,
Lemma \ref{L17} implies that if the subgroups of $G$ belonging to $[Stab_G(x),G]$
commute with each other for some $x\in X$ then the subgroups of $G$ belonging to
$[Stab_G(y),G]$ commute with each other for all $y\in X$.
\end{remark}

\begin{construction}

Let $G$ be a group and $G_a$ be a subgroup of $G$ such that $HK=KH$ is
satisfied for all subgroups $H,K$ of $G$ containing $G_a$. Let
$N^*$ denote the right quotient set $G/G_a$, that is, the set of all right cosets
$G_ag$ $(g\in G$) of $G$ defined by $G_a$. Let $S=G\cup N^*\cup \{ 0\}$, where
$0$ is a symbol not contained in $G\cup N^*$. On $S$ we define an operation as follows.
If $g,h\in G$ then let $gh$ be the original product of $g$ and $h$ in $G$. If
$a\in N$ then, for arbitrary $s\in S$, let $sa=0$. For arbitrary $g\in G$ and arbitrary $G_ah\in N^*$, let
$(G_ah)g=G_a(hg)$. It is easy to check that $S$ is a semigroup.
\end{construction}

\begin{theorem}\label{T2} A finite semigroup $S$ is a permutable semigroup which is a
semilattice of a group $G$ and a nil semigroup such that
the identity element of $G$ is a right [left] identity element of $S$ and
$SN=\{ 0\}$ [$NS=\{ 0\}$] if and only if it is isomorphic to a semigroup defined in
Construction 1 [the dual of Construction 1].
\end{theorem}

\noindent
{\bf Proof}. First of all we show that the semigroup $S$ defined in Construction 1, is
a permutable semigroup. It is clear that $S$ is
a semilattice of the group $G$ and the null semigroup $N=N^*\cup \{ 0\}$
such that $SN=\{ 0\}$ and the identity element $e$ of $G$ is a right
identity element of $S$. Moreover, $(G_ag)G=N^*$ for all $G_ag\in N^*$.
Thus $N^*$ is a right $G$-set and $G$ acts on $N^*$ tranisitively. By Lemma 4.20 of
\cite{McKenzie}, the congruence lattice
$\it{Con}(N^*)$ of the $G$-set $N^*$ is isomorphic to the interval $[\it{Stab_G}(G_a),G]$,
where $\it{Stab_G}(G_a)=\{g\in G; G_ag=G_a\}=G_a$.
Let $\alpha$ be a non-universal congruence of $S$. Then, by Lemma \ref{L16}, $[g]_{\alpha}\subseteq G$
for every $g\in G$ and $[0]_{\alpha}$ is either $\{ 0\}$ or $N$.
As $N^*$ is a right $G$-set and $G$ acts on $N^*$ transitively, moreover
the restriction $\alpha ^*$ of $\alpha$ to $N^*$ is in $Con(N^*)$, there is a subgroup
$H_{\alpha ^*} \in [G_a,G]$ which determines $\alpha$ on $N^*$.

Let $\alpha$ and $\beta$ be arbitrary congruences of $S$. We show that
$\alpha \circ \beta =\beta \circ \alpha$. We can suppose that $\alpha$ and $\beta$
are not the universal relations of $S$.
Assume $(b,c)\in \alpha \circ \beta$ for arbitrary elements $b$ and $c$ of $S$. Then there
is an element $x\in S$ such that $(b,x)\in \alpha$, $(x,c)\in \beta$.
We have two cases.

\noindent
Case 1: $x\in G$. Then, by Lemma \ref{L16}, $b,c\in G$. As every group is permutable,
there is an element $y\in G$ such that $(b,y)\in \beta$ and $(y,c)\in \alpha$.
Hence $(b,c)\in \beta \circ \alpha$.

\noindent
Case 2: $x\in N$. Then, by Lemma \ref{L16}, $b,c\in N$.
If $[0]_{\alpha}=N$ or $[0]_{\beta}=N$ then $(b,c)\in \alpha$ or
$(b,c)\in \beta$ and so $(b,c)\in \beta \circ \alpha$.
Consider the case when $[0]_{\alpha}=[0]_{\beta}=\{ 0\}$. Then $N^*$ is saturated
by both $\alpha$ and $\beta$. If $x=0$ then $b=c=0$ and so
$(b,c)\in \beta \circ \alpha$. If $x\in N^*$ then $b,c\in N^*$. If
$\alpha ^*$ and $\beta ^*$ denote the restriction of $\alpha$ and $\beta$ to $N^*$,
respectively, then $H_{\alpha ^*}, H_{\beta ^*}\supseteq G_a$. As
$H_{\alpha ^*}H_{\beta ^*}=H_{\beta ^*}H_{\alpha ^*}$, we get
$\alpha ^* \circ \beta ^* =\beta ^*\circ \alpha ^*$ by Lemma \ref{L17}. Hence
$(b,c)\in \beta \circ \alpha$.

Thus we have $(b,c)\in \beta \circ \alpha$ in both cases.
Consequently, $\alpha \circ \beta \subseteq \beta \circ \alpha$. By the symmetry,
we get $\alpha \circ \beta =\beta \circ \alpha$. Thus $S$ is a permutable semigroup.

Conversely, assume that $S$ is a permutable semigroup which is a semilattice
of a group $G$ and a non-trivial nil semigroup $N$ such that the identity element of $G$
is a right identity element of $S$ and $SN=\{ 0\}$.
Then $N$ is a null semigroup and $aG=N^*$ for every $a\in N^*$
by Lemma \ref{L15}. Thus $N^*$ is a right $G$-set and $G$ acts on $N^*$ transitively.
Fix an element $a$ in $N^*$ and consider $G_a=Stab_G(a)=\{ g\in G:\ ag=a\}$.
It is easy to check that $ag=ah$ for some $g,h\in G$ if and only if $G_ag=G_ah$.
Thus $|N^*|=|G:G_a|$. Let $\phi$ be the bijection of $N^*$ to the factor set $G/G_a$
defined by $\phi (b)=G_ag$ if $b=ag$. It is clear that $\phi$ is well defined.
Moreover, for all $g,h\in G$, $(G_ag)h=G_a(gh)$ implies $(\phi (b))h=\phi (bh)$.
If we identify every $b\in N^*$ with $\phi (b)$ then $N^*$ can be considered as the set
of all right cosets of $G$ defined by $G_a$, and the operation on $S$ is defined
as in the Construction 1.
Let $H$ and $K$ be arbitrary subgroups of $G$ containing the subgroup $G_a$.
Let $\alpha '_H=\alpha _H \cup 1_S$ and $\alpha '_K=\alpha _K \cup 1_S$,
where $\alpha _H=\psi (H)$ and $\alpha _K=\psi (K)$ are congruences of the right
$G$-set $N^*$ defined by $H$ and $K$,
respectively (for $\psi$, we refer to Remark \ref{R4}). It is easy to see that $\alpha '_H$
and $\alpha '_K$ are congruences of $S$. As $S$ is permutable,
they commute with each other from which we get $\alpha _H\circ \alpha _K=
\alpha _K\circ \alpha _H$. Hence $HK=KH$ by Lemma \ref{L17}. Thus the theorem is proved.
\hfill\openbox

\subsection{When the identity element of $G$ is the two-sided identity element of $S$}

\begin{lemma}\label{L18} If $S$ is a permutable semigroup which is a semilattice of a group $G$
and a non-trivial nilpotent semigroup $N$ of nilpotency degree $t$ such that the identity element of $G$ is an identity element of $S$
then, for all $a\in N^k-N^{k+1}$, $GaG=N^k-N^{k+1}$ is satisfied for every $k=1, \dots , t-1$.
\end{lemma}

\noindent
{\bf Proof}.
Let $a\in N^k-N^{k+1}$ be arbitrary. It is clear that $GaG\subseteq N^k-N^{k+1}$.
As the ideals of $S$ form a chain, $N^{k+1}\subseteq SaS$. It is clear that $SaS\subseteq N^k$.
Assume $SaS\neq N^k$ for every $a\in N^k-N^{k+1}$. As $S$ is finite, there is an
element $b\in N^k-N^{k+1}$ such that $SaS\subseteq SbS\neq N^k$,
Let $c\in N^k-SbS$ be arbitrary. Then $SbS\subset ScS$ which is a contradiction.
Consequently, $SaS=N^k$ for some $a\in N^k-N^{k+1}$. Thus
$N^k=SaS=GaG\cup GaN\cup NaG\cup NaN\subseteq GaG\cup N^{k+1}$ which implies
$GaG=N^k-N^{k+1}$.\hfill\openbox

\begin{lemma}\label{L19} Let the finite semigroup $S$ be a semilattice of a group $G$
and a non-trivial nilpotent semigroup $N$ of nilpotency degree $t$ such that,
for every $a\in N^i-N^{i+1}$ ($i=1, \dots t-1$), $GaG=N^i-N^{i+1}$ is satisfied.
Then the ideals of $S$ are $S, N, N^2, \dots N^t=\{ 0\}$.
\end{lemma}

\noindent
{\bf Proof}. It is clear that $S, N, N^2, \dots , N^t=\{ 0\}$ are ideals of $S$.
Let $I$ be an arbitrary ideal of $S$. Let $j$ be the least positive integer such that
$I\cap N^j\neq \emptyset$. If $a\in I\cap (N^j-N^{j+1})$ then
$N^j-N^{j+1}=GaG\subseteq I$. Let $b\in N^{j+1}-N^{j+2}$ supposing that $N^{j+1}\neq \{ 0\}$.
There are elements $x_1,\dots , x_{j+1}\in N-N^2$ such that
$b=x_1\dots x_{j+1}$. It is clear that $x_1\dots x_j\in N^j-N^{j+1}\subseteq I$
and so $b\in I$ which implies that $N^{j+1}-N^{j+2}\subseteq I$. Continouing this procedure,
we get that $N^j=I\cap N$. If $I\cap G=\emptyset$ then $I=N^j$. Assume that
$I\cap G\neq \emptyset$. Then $G\subseteq I$. Moreover, for all $i=1,\dots ,t-1$,
and every $a\in N^i-N^{i+1}$, $N^i-N^{i+1}=GaG\subseteq I$ which implies that $I=S$.
Thus the lemma is proved.\hfill\openbox

\begin{lemma}\label{L20} Let $S$ be a semigroup which is a semilattice of a group $G$
and a nilpotent semigroup $N$ of nilpotency degree $t$ such that,
for every $i\in \{1, \dots , t-1\}$ and for some (and so for every) $a\in N^i-N^{i+1}$,
$GaG=N^i-N^{i+1}$ is satisfied.
Then, for every non-universal congruence $\alpha$ of $S$, $[0]_{\alpha}=N^j$
for some positive integer $j=1, \dots ,t$
and $[g]_{\alpha}\subseteq G$ for every $g\in G$, moreover
$[a]_{\alpha}\subseteq N^i-N^{i+1}$ for every $a\in N^i-N^{i+1}$ ($i=1, \dots j-1$).
\end{lemma}

\noindent
{\bf Proof}. Let $\alpha$ be a non-universal congruence of $S$. If $(g,a)\in \alpha$ for some $g\in G$ and $a\in N$
then $(g^t,a^t)\in \alpha$. As $a^t=0$, and $(ugv,0)\in \alpha$ for all $u,v\in G$,
we get $G\subseteq [0]_{\alpha}$. Let $a\in N$ be an arbitrary element.
Then $(gah,0)\in \alpha$ for all $g,h\in G$ and so $N^i-N^{i+1}\subseteq [0]_{\alpha}$
for everi $i=1, \dots t-1$. Thus $S=[0]_{\alpha}$ which is a contradiction.
Hence $[g]_{\alpha}\subseteq G$ for every $g\in G$.
By Lemma \ref{L19}, the ideals of $S$ are
$S, N, N^2, \dots N^{t-1}, N^t=\{ 0\}$. Then there is a least positive integer
$j\in \{1,2, \dots t\}$ such that $[0]_{\alpha}=N^j$. If $j=1$ or $j=2$ then the assertion
is true for $\alpha$. Assume $j\geq 3$.
Let $a\in N^{j-1}-N^j$ be arbitrary. It is clear that $(a,b)\notin \alpha$ for every
$b\in N^j$. Assume $(a,b)\in \alpha$ for some $b\in N^{k-1}-N^k$ for some
$k<j$.
There are elements $x_1, \dots ,x_{j-1}\in N-N^2$ such that
$a=x_1\dots x_{k-1}\dots x_{j-1}$. It is clear that $x_1\dots x_{j-2}\in N^{j-2}$,
$x_1\dots x_{j-3}\in N^{j-3}-N^{j-2}$, and finally, $c=x_1\dots x_k-1\in N^{k-1}-N_k=GbG$.
Then $c=gbh$ for some $g,h\in G$.
Thus $(c,gah)\in \alpha$. As $gah\in N^{j-1}-N^j$, $d=gahx{k}\dots x_{j-1}\in N^j$
and so $a=cx_{k+1}\dots x_{j-1}$ implies $(a,d)\in \alpha$ which is impossible.
Hence $[a]_{\alpha}\subseteq N^{j-1}-N^j$. Thus the lemma is proved.\hfill\openbox

\medskip

For an arbitrary group $G$, let $G^*$ denote the dual of $G$, that is, $xy=u$
in $G^*$ if and only if $yx=u$ in $G$.

\begin{theorem}\label{T3}
Let $S$ be a finite semigroup which is a semilattice of a group $G$
and a non-trivial nilpotent semigroup $N$ of nilpotency degree $t$ such that
the identity element of $G$ is the identity element of $S$.
Then $S$ is permutable if and only if,
for all $i=1, \dots , t-1$,
there is an element $a_i$ in $N^i-N^{i+1}$ such that
$Ga_iG=N^i-N^{i+1}$, and $HK=KH$ is satisfied
for all subgroups $H,K\supseteq G_{a_i}=\{ (g,h)\in G^*\times G:\ ga_ih=a_i\}$.
\end{theorem}

\noindent
{\bf Proof}. Let $S$ be a finite semigroup which is a semilattice of a
group $G$ and a nilpotent semigroup $N$ of nilpotency degree $t$ such that
the identity element of $G$ is the identity element of $S$. First assume that
$S$ is permutable.
Let $i\in \{ 1, \dots ,t-1\}$ be arbitrary. Then, for every $a_i\in N^i-N^{i+1}$,
$Ga_iG=N^i-N^{i+1}$ is satisfied by Lemma \ref{L18}. It is a matter of checking
to see that this result implies that $N^i-N^{i+1}$ is a right $(G^*\times G)$-set
($a(g,h)=gah$ for every $a\in N^*$ and every $(g,h)\in G^*\times G$)
and $G^*\times G$ acts on $N^i-N^{i+1}$ transitively.
Let $G_{a_i}=Stab_{G^*\times G}(a_i)=\{(g,h)\in G^*\times G:\ ga_ih=a_i\}$. By Lemma 4.20
of \cite{McKenzie} the congruence lattice $Con(N^i-N^{i+1})$ of the right
$(G^*\times G)$-set $N^i-N^{i+1}$ is isomorphic to $[Stab_{G^*\times G}(a_i),G^*\times G]$.
The corresponding
isomorphisms $\phi:\ \alpha _i\ \mapsto \ H_{\alpha _i}$ ($\alpha _i\in Con(N^i-N^{i+1}$) and
$\psi :\ H \mapsto \alpha ^{(i)}_H$ ($H\in [Stab_{G^*\times G}(a_i),G^*\times G]$)
defined as in Remark \ref{R4}.
Let $H$ be an arbitrary subgroup of
$G^*\times G$ containing the subgroup $G_{a_i}$. Let $\alpha '_H$ be the relation of
$S$ defined by $(a,b)\in \alpha '_H$ if and only if $a=b$ or $a,b\in N^{i+1}$ or
$a,b\in N^i-N^{i+1}$ and $(a,b)\in \alpha ^{(i)}_H$.
It is clear that $\alpha '_H$ is an equivalence relation.
We show that it is a congruence of $S$.
Assume $(a,b)\in \alpha '_H$ for some $a,b\in S$. We can suppose that $a\neq b$.
If $a,b\in N^{i+1}$ then $sa, sb, as, bs\in N^{n+1}$ and so $(sa,sb)\in \alpha '_H$
and $(as,bs)\in \alpha '_H$. Consider the case when $a, b\in N^i-N^{i+1}$.
Then $(a,b)\in \alpha ^{(i)}_H$ and so, for every $x\in G$, we have
$(a(e,x),b(e,x))\in \alpha ^{(i)}_H$ and $(a(x,e),b(x,e))\in \alpha ^{(i)}_H$,
because $\alpha ^{(i)}_H$ is a congruence of the right $(G^*\times G)$-set
$N^i-N^{i+1}$. Thus $(ax,bx)\in \alpha ^{(i)}_H$ and $(xa,xb)\in \alpha ^{(i)}_H$. Hence
$(ax,bx)\in \alpha '_H$ and $(xa,xb)\in \alpha '_H$.
If $u\in N$ then
$ua,ub, au, bu \in N^{i+1}$ and so $(au,bu)\in \alpha '_H$ and $(ua,ub)\in \alpha '_H$.
Consequently, $\alpha '_H$ is a congruence on $S$.
Let $H$ and $K$ be arbitrary subgroups of $G^*\times G$ containing the subgroup
$G_{a_i}$. Let $\alpha '_H$ and $\alpha '_K$ be the congruences of $S$ defined by $H$ and $K$
(see above). As $S$ is permutable,
$\alpha '_H\circ \alpha '_K=\alpha '_K\circ \alpha '_H$ from
which we get $\alpha ^{(i)}_H\circ \alpha ^{(i)}_K=\alpha ^{(i)}_K\circ \alpha ^{(i)}_H$. Then
$HK=KH$ by Lemma \ref{L17}. Thus the necessity of the permutability of $S$ is proved.

Conversely, assume that,
for all $i=1, \dots , t-1$,
there is an element $a_i$ in $N^i-N^{i+1}$ such that
$Ga_iG=N^i-N^{i+1}$, and $HK=KH$ is satisfied
for all subgroups $H,K\supseteq G_{a_i}=\{ (g,h)\in G^*\times G:\ ga_ih=a_i\}$.
We note that, from $Ga_iG=N^i-N^{i+1}$, it follows that $GaG=N^i-N^{i+1}$ for every
$a\in N^i-N^{i+1}$. Thus $N^i-N^{i+1}$ is a right $(G^*\times G)$-set and
$G^*\times G$ acts on $N^i-N^{i+1}$ transitively. By Lemma 4.20 of \cite{McKenzie},
the congruence lattice $Con(N^i-N^{i+1})$ of the $(G^*\times G)$-set $N^i-N^{i+1}$ is
isomorphic to $[Stab_{G^*\times G}(a_i),G^*\times G]$ (for the corresponding
isomorphisms we refer to Remark \ref{R4}). By Lemma \ref{L19},
the ideals of $S$ are $S, N, N^2, \dots ,N^t$.
Let $\alpha$ be a non-universal congruence on $S$. Then, by Lemma \ref{L20},
$[0]_{\alpha}=N^j$ for some positive integer $j\in \{ 1, \dots t\}$,
$[g]_{\alpha}\subseteq G$ for every $g\in G$, and
$[a]_{\alpha}\subseteq N^i-N^{1+1}$ for every $a\in N^i-N^{i+1}$ ($i=1, \dots j-1$).
Let $\alpha _i$ denote the restriction of $\alpha$ to $N^i-N^{i+1}$, and let
$H^{(i)}_{\alpha}=\phi (\alpha _i)=\{ (g,h)\in G^*\times G:\ a_i(g,h)\ \alpha _i \ a_i\}$
($a_i\in N^i-N^{i+1},\ i=1,\dots t-1$).
$H^{(i)}_{\alpha}$ is a subgroup of $G^*\times G$ and $G_{a_i}\subseteq H^{(i)}_{\alpha}$.
Let $\beta$ be an arbitrary non-universal congruence on $S$. As $G_{a_i}\subseteq H^{(i)}_{\beta}$,
we have $H^{(i)}_{\alpha}H^{(i)}_{\beta}=H^{(i)}_{\beta}H^{(i)}_{\alpha}$.
As $\alpha _i$ and $\beta _i$ are in the congruence lattice
$Con(N^i-N^{i+1})$ of the right $(G^*\times G)$-set $N^i-N^{i+1}$, we have
$\alpha _i\circ \beta _i=\beta _i\circ \alpha _i$. We show that $\alpha \circ \beta =\beta \circ \alpha$.
Assume $(a,b)\in \alpha \circ \beta$ for some $a,b\in S$.
Then there is an element $c\in S$ such that $(a,c)\in \alpha$ and $(c,b)\in \beta$.
If $c\in G$ then $a,b\in G$ by Lemma \ref{L20}. As every group is permutable,
we get
$(a,b)\in \beta \circ \alpha$. By Lemma \ref{L2}, $[0]_{\alpha}\subseteq [0]_{\beta}$ or
$[0]_{\beta}\subseteq [0]_{\alpha}$. Assume $[0]_{\alpha}\subseteq [0]_{\beta}$.
If $c\in [0]_{\beta}$ then $a,b\in [0]_{\beta}$ and so $(a,b)\in \beta$,
$(b,b)\in \alpha$ implies $(a,b)\in \beta \circ \alpha$.
Assume $c\notin [0]_{\beta}$, $c\in N^i-N^{i+1}$. Then, by Lemma \ref{L20},
$a,b\in N^i-N^{i+1}$. Thus $(a,b)\in \alpha _i\circ \beta _i=\beta _i\circ \alpha _i$ (see
also Lemma \ref{L17}).
Thus $(a,b)\in \beta \circ \alpha$.
The proof of $(a,b)\in \beta \circ \alpha$ is similar in that case when
$[0]_{\beta}\subseteq [0]_{\alpha}$. Thus
$\alpha \circ \beta \subseteq \beta \circ \alpha$.
The proof of $\beta \circ \alpha \subseteq \alpha \circ \beta$ is similar.
Thus $\alpha \circ \beta =\beta \circ \alpha$. Hence $S$ is permutable.
\hfill\openbox

\medskip

By Remark \ref{R5}, if a semigroup satisfies the conditions of Theorem \ref{T3}
then $Ga_iG=N^i-N^{i+1}$ and $HK=KH$ is satisfied for all subgroups $H,K\supseteq G_{a_i}=
\{ (g,h)\in (G^*\times G):\ ga_ih=a_i\}$, for all elements $a_i$ in $N^i-N^{i+1}$.

\medskip

The authors wish to thank the referee for his valuable suggestions for improving
the first version of the paper, in particular, for finding much shorter proofs of
Theorem 1 and Lemma 14.

\bigskip

\noindent
A. De\'ak and A. Nagy,

\noindent
Department of Algebra,

\noindent
Institute of Mathematics,

\noindent
Budapest University of Technology and Economics,

\noindent H-1111 Budapest, M\H{u}egyetem rkp. 9.,Hungary;

\noindent
{\it e-mail}: A. De\'ak: ignotus@math.bme.hu, A. Nagy: nagyat@math.bme.hu

\end{document}